\documentclass[11pt]{amsart}
\usepackage{amsmath,amssymb,mathtools}
\usepackage[margin=1.15in]{geometry}
\usepackage[hidelinks]{hyperref}
\newtheorem{theorem}{Theorem}[section]
\newtheorem{lemma}[theorem]{Lemma}
\newtheorem{corollary}[theorem]{Corollary}
\theoremstyle{remark}

\newcommand{\R}{\mathbb R}
\newcommand{\Z}{\mathbb Z}
\newcommand{\dimH}{\dim_{\mathrm H}}
\newcommand{\dimP}{\dim_{\mathrm P}}
\numberwithin{equation}{section}
\title[Dimensions of successive-minimum windows]{Exact dimensions for $k$-singularity and arbitrary windows of successive minima}
\author{Shahriar Mirzadeh}
\address{Department of Mathematics, Brandeis University,
415 South Street, Waltham, MA 02453, USA}
\email{shahmir@brandeis.edu}
\date{November 2026}
\begin{document}
\begin{abstract}
We determine the exact Hausdorff and packing dimensions of
sets of matrices for which two prescribed successive minima
of the associated diagonal lattice orbit tend to zero and
infinity, respectively. When the prescribed minima are
separated by exactly one intermediate minimum, our formula
resolves a conjecture of Das, Fishman, Simmons and
Urba\'nski on $k$-singular matrices.

More generally, we prescribe the lower asymptotic density of
times at which the two minima are simultaneously below
$\varepsilon$ and above $\varepsilon^{-1}$, respectively, in the
limit as $\varepsilon\to0$. We determine the complete Hausdorff
and packing dimension spectra and show that both are affine
functions of the prescribed density. The same formulas hold for
level sets and superlevel sets, even after excluding trivially
singular matrices. At density one, the dimensions agree with
those for permanent divergence except when the two prescribed
minima are adjacent, in which case they are strictly larger.

The principal upper bounds follow from the variational principle
in parametric geometry of numbers, using a piecewise-linear
potential and a block-merging argument that reduces general
templates to two-block configurations. Explicit two-block
excursions give the matching lower bounds when the divergence
proportion is prescribed. Permanent divergence between adjacent
successive minima is treated separately using convexity of
partial sums.
\end{abstract}
\maketitle

\section{Introduction}
Fix integers $m,n\geq1$ and put $d=m+n\geq3$. For $A\in M_{m,n}(\R)$ and $t\geq0$, let
\[
 \Lambda_A(t)=g_tu_A\Z^d,\qquad
 g_t=\operatorname{diag}(e^{t/m}I_m,e^{-t/n}I_n),\qquad
 u_A=\begin{pmatrix}I_m&A\\0&I_n\end{pmatrix}.
\]
Since $\det(g_tu_A)=1$, each $\Lambda_A(t)$ is unimodular. Write $\lambda_i(\Lambda)$ for the $i$th successive minimum of a lattice $\Lambda$ with respect to the sup-norm unit ball. The vector with integer label $(p,q)\in\Z^m\times\Z^n$ has coordinates $(e^{t/m}(p+Aq),e^{-t/n}q)$. Thus $\lambda_i(\Lambda_A(t))$ is the smallest scale $\varepsilon$ at which there are $i$ linearly independent integer pairs satisfying $\|p+Aq\|_\infty\leq\varepsilon e^{-t/m}$ and $\|q\|_\infty\leq\varepsilon e^{t/n}$. Put
\[
 h_A(t)=(h_{A,1}(t),\ldots,h_{A,d}(t)),\qquad
 h_{A,i}(t)=\log\lambda_i(\Lambda_A(t)).
\]
All limiting conditions below are independent of the chosen norm, since equivalent norms change $h_A$ by a bounded amount.

By the Dani correspondence \cite{Dani}, $A$ is singular if and only if
$h_{A,1}(t)\to-\infty$; see also \cite[Theorem~3.2]{DFSU}.
Cheung \cite{Cheung} determined the Hausdorff dimension of the set of
singular pairs, and Cheung and Chevallier \cite{CC} extended this result
to singular vectors in arbitrary dimension greater than one.
Kadyrov, Kleinbock, Lindenstrauss and Margulis \cite{KKLM} obtained
the sharp Hausdorff dimension upper bound for singular matrices,
while Das, Fishman, Simmons and Urba\'nski \cite{DFSU} established
the matching lower bound and the corresponding packing dimension formula.

The study of divergence conditions involving higher
successive minima goes back to a conjecture of Schmidt,
proved by Moshchevitin \cite{Moshchevitin}.

Das, Fishman, Simmons and Urba\'nski \cite[Definition~3.13]{DFSU} call $A$ \emph{$k$-singular} if $h_{A,k-1}(t)\to-\infty$ and $h_{A,k+1}(t)\to+\infty$, where $2\leq k\leq d-1$. They announced an exact dimension formula for this set in~\cite[Theorem~1.12]{DFSU2017}. In the subsequent full paper \cite[Theorem~3.14 and Conjecture~5.1]{DFSU}, the authors established the Hausdorff dimension lower bound and conjectured the matching upper bound for both Hausdorff and packing dimension. The first version of the present work, posted as \href{https://arxiv.org/abs/2611.30847v1}{arXiv:2611.30847v1} on November 23, 2026, proved the matching packing-dimension upper bound by a block-merging argument. Here we extend that argument to arbitrary windows of indices and to divergence on average, using a simpler choice of potential.

We consider a more general form of divergence in which
two prescribed successive minima are required to tend
to zero and infinity, respectively. For integers
$1\leq a\leq b\leq d-1$, define
\begin{equation}\label{eq:window}
 \mathcal S_{a,b}
 =
 \left\{
 A\in M_{m,n}(\mathbb R):
 h_{A,a}(t)\to-\infty,\quad
 h_{A,b+1}(t)\to+\infty
 \right\},
\end{equation}
where both limits are taken as $t\to\infty$.

We say that the orbit associated with $A$ exhibits
\emph{permanent divergence in the window $[a,b]$}
if $A\in \mathcal S_{a,b}$. Equivalently, for every $R>0$,
\[
 h_{A,a}(t)\leq-R,\qquad h_{A,b+1}(t)\geq R
\]
for all sufficiently large $t$.

By the ordering of the successive minima, these conditions
force the first $a$ logarithmic minima to tend to
$-\infty$ and the last $d-b$ logarithmic minima to tend
to $+\infty$. Consequently, for all sufficiently large
$t$, the number of strictly negative logarithmic minima
belongs to $\{a,\ldots,b\}$. This motivates calling the
integer interval $[a,b]\cap\mathbb Z$ the \emph{window
of indices}.

We refer to $[1,d-1]$ as the \emph{full window} and
to a window of the form $[j,j]$ as a \emph{singleton
window}. The latter corresponds to prescribing the
divergence of two adjacent successive minima.

Following \cite{DFSU}, we call a matrix $A$ trivially
singular if
\[
 h_{A,i+1}(t)-h_{A,i}(t)\to+\infty
\]
for some fixed $1\leq i<d$. For any set of matrices
$\mathcal E$, we write $\mathcal E^*$ for the subset
obtained by removing all trivially singular matrices.

With $\{x\}=x-\lfloor x\rfloor$, set
\begin{equation}\label{eq:F}
 F_{m,n}(j)=mn-\frac{mn\,j(d-j)}{d^2}
 -\left\{\frac{jm}{d}\right\}\left\{\frac{jn}{d}\right\},
 \qquad 1\leq j\leq d-1,
\end{equation}
and write $D_{a,b}=\max\{F_{m,n}(a),F_{m,n}(b)\}$. The function in \eqref{eq:F} is the one appearing in the lower bound of \cite[Theorem~3.14]{DFSU}.

\begin{theorem}\label{thm:main}
For $1\leq a<b\leq d-1$, both $\mathcal S_{a,b}$ and $\mathcal S_{a,b}^*$ have Hausdorff and packing dimension
\begin{equation}\label{eq:main}
 \max_{a\leq j\leq b}F_{m,n}(j)=D_{a,b},
\end{equation}
where the maximum on the left is over integer indices.
\end{theorem}
Several important cases are contained in Theorem~\ref{thm:main}.
For the window $[a,b]=[k-1,k]$, the conditions defining
$\mathcal S_{a,b}$ become
\[
 h_{A,k-1}(t)\to-\infty,\qquad
 h_{A,k+1}(t)\to+\infty.
\]
Thus $\mathcal S_{k-1,k}$ is precisely the set of $k$-singular
matrices, and Theorem~\ref{thm:main} recovers the dimension formula
established in the first version of this paper.

For the full window $[1,d-1]$, the defining conditions are
$h_{A,1}(t)\to-\infty$ and $h_{A,d}(t)\to+\infty$.
By Minkowski's second theorem, the first condition implies
the second. Hence $\mathcal S_{1,d-1}$ is precisely the set of
singular matrices. In this case,
\[
 F_{m,n}(1)=F_{m,n}(d-1)=mn-\frac{mn}{d},
\]
and Theorem~\ref{thm:main} recovers the exact dimension formula for
singular matrices established in \cite[Theorem~3.1]{DFSU}.

Singleton windows exhibit a different behavior.
For $a=b=j$, the defining conditions become
\[
 h_{A,j}(t)\to-\infty,\qquad
 h_{A,j+1}(t)\to+\infty.
\]
In particular, the gap between the $j$th and $(j+1)$st
logarithmic minima tends to infinity, so every matrix
in $\mathcal S_{j,j}$ is trivially singular. The dimension of
this set is given by the following theorem.

\begin{theorem}\label{thm:singleton}
For $1\leq j\leq d-1$, put $p_j=\lceil mj/d\rceil-1$. Then
\begin{equation}\label{eq:singleton}
 \dimH\mathcal S_{j,j}=\dimP\mathcal S_{j,j}=mn-(j-p_j)(m-p_j).
\end{equation}
Moreover, $\mathcal S_{j,j}^*=\varnothing$.
\end{theorem}

We next consider a time-average version of the preceding
divergence conditions. Rather than requiring the two
prescribed minima to remain beyond every fixed threshold
for all sufficiently large times, we prescribe the
lower asymptotic proportion of time for which both conditions
hold simultaneously, and then let the thresholds tend to infinity.

For a continuous map $h:[0,\infty)\to\mathbb R^d$,
a threshold $R>0$, and $T>0$, define
\[
 P_{a,b,R}(h,T)
 =
 \frac{1}{T}
 \left|
 \left\{
 t\in[0,T]:
 h_a(t)\leq-R,\quad h_{b+1}(t)\geq R
 \right\}
 \right|,
\]
where $|\cdot|$ denotes Lebesgue measure. Thus
$P_{a,b,R}(h,T)$ is the proportion of time in $[0,T]$
during which the two prescribed logarithmic minima lie
beyond the thresholds $-R$ and $R$, respectively.

Define the \emph{divergence proportion in the window
$[a,b]$} by
\[
 P_{a,b}(h)
 =
 \lim_{R\to\infty}
 \liminf_{T\to\infty}P_{a,b,R}(h,T).
\]
For each fixed threshold $R$, the inner limit inferior
measures the lower asymptotic density of such times.
These densities are nonincreasing in $R$ and take
values in $[0,1]$, so the outer limit exists.

We say that the orbit associated with $A$ exhibits
\emph{divergence on average in the window $[a,b]$}
if
\[
 P_{a,b}(h_A)=1.
\]
Equivalently, for every $R>0$, the proportion of time
during which
\[
 h_{A,a}(t)\leq-R,\qquad h_{A,b+1}(t)\geq R
\]
tends to one as $T\to\infty$.

Permanent divergence in a window implies divergence
on average in that window, but the converse need not
hold: divergence on average allows returns to regions
where the prescribed inequalities fail, provided that
such returns occupy an asymptotically negligible
proportion of time.

More generally, for $\rho\in[0,1]$, define the level
and superlevel sets
\[
 \mathcal L_{a,b}(\rho)
 =
 \{A\in M_{m,n}(\mathbb R):P_{a,b}(h_A)=\rho\},
\]
and
\[
 \mathcal U_{a,b}(\rho)
 =
 \{A\in M_{m,n}(\mathbb R):P_{a,b}(h_A)\geq\rho\}.
\]
In particular,
\[
 \mathcal S_{a,b}\subseteq
 \mathcal L_{a,b}(1)=\mathcal U_{a,b}(1).
\]
For the full window $[1,d-1]$, Minkowski's second theorem
implies that
\[
 P_{1,d-1}(h_A)
 =\lim_{R\to\infty}\liminf_{T\to\infty}
 \frac1T\big|\{t\in[0,T]:h_{A,1}(t)\leq-R\}\big|.
\]
Indeed, there is a constant $C_d\geq0$ such that
$|\sum_{i=1}^d h_{A,i}(t)|\leq C_d$ for every $t$. By the
ordering of the logarithmic minima,
\[
 h_{A,d}(t)\geq\frac{-h_{A,1}(t)-C_d}{d-1}.
\]
Thus $h_{A,1}(t)\leq-(d-1)R-C_d$ implies both
$h_{A,1}(t)\leq-R$ and $h_{A,d}(t)\geq R$. Conversely,
the simultaneous inequalities imply $h_{A,1}(t)\leq-R$.
Taking the lower asymptotic proportions and then letting
$R\to\infty$ gives the stated identity.
The quantity on the right is the proportion of time near
infinity defined in \cite[Section~3.4]{DFSU}.

The following theorem determines the complete Hausdorff
and packing dimension spectra of these sets.
\begin{theorem}\label{thm:spectrum}
For $1\leq a\leq b\leq d-1$ and $0\leq\rho\leq1$, the sets $\mathcal L_{a,b}(\rho)$, $\mathcal U_{a,b}(\rho)$ and their starred subsets all have Hausdorff and packing dimension
\begin{equation}\label{eq:spectrum}
 (1-\rho)mn+\rho D_{a,b}.
\end{equation}
\end{theorem}

For the full window $[1,d-1]$, Theorem~\ref{thm:spectrum}
recovers the dimension formula of \cite[Theorem~3.11]{DFSU}
and extends it to arbitrary windows of successive minima.

\medskip
Theorems~\ref{thm:main}--\ref{thm:spectrum} reveal a
distinction between singleton and nonsingleton windows.
For $a<b$, permanent divergence and divergence on average
give rise to sets of the same Hausdorff and packing dimension.
For singleton windows, however, both dimensions increase
strictly when permanent divergence is replaced by divergence
on average. Indeed, writing $\theta=mj/d-p_j\in(0,1]$,
Theorems~\ref{thm:singleton} and~\ref{thm:spectrum} give
\[
 \dimH\mathcal L_{j,j}(1)-\dimH\mathcal S_{j,j}
 =\theta(m+j-2p_j-1)>0,
\]
where the positivity follows from
$p_j+1\leq\min\{m,j\}$.

The distinction has a simple template-level explanation. For
permanent divergence in a singleton window, the partial sum
$H_j$ is eventually convex and tends to $-\infty$; because its
allowed slopes form the finite set $Z(j)$, its slope is eventually
constant and negative. In the divergence-on-average construction,
returns to coincident coordinates allow this slope to reset between
excursions, producing the larger average contraction rate.

For example, when $(m,n)=(1,5)$ and $j=3$, we have
\[
 \dimH\mathcal S_{3,3}=2,\qquad
 \dimH\mathcal L_{3,3}(1)
 =\dimH\mathcal L_{3,3}^*(1)=\frac72,
\]
whereas $\mathcal S_{3,3}^*=\varnothing$.

\medskip
\noindent\textbf{Higher-order singularity.}
Xing \cite[Definition~1.14]{Xing} introduced singular matrices of
order $r$, which in our notation form the set
\[
 \operatorname{Sing}_r(m,n)
 =\{A\in M_{m,n}(\R):h_{A,r}(t)\to-\infty\}.
\]
The computation of its exact Hausdorff and packing dimensions
was left open in \cite[Remark~1.19]{Xing}. By the ordering of the
logarithmic minima and Minkowski's second theorem,
$\operatorname{Sing}_r(m,n)=\mathcal S_{r,d-1}$ for
$1\leq r\leq d-1$. Indeed, if $h_{A,r}(t)\to-\infty$, then
\[
 (d-r)h_{A,d}(t)\geq-rh_{A,r}(t)-O_{m,n}(1),
\]
so $h_{A,d}(t)\to+\infty$. Theorems~\ref{thm:main} and
\ref{thm:singleton} therefore give
\[
 \dimH\operatorname{Sing}_r(m,n)
 =\dimP\operatorname{Sing}_r(m,n)
 =\begin{cases}
 mn-\dfrac{mn}{d},&1\leq r\leq d-2,\\[4pt]
 mn-n,&r=d-1.
 \end{cases}
\]
Here the first case also follows directly from
\cite[Theorems~3.1 and~3.14]{DFSU}, since
$\mathcal S_{d-2,d-1}\subseteq\operatorname{Sing}_r(m,n)
\subseteq\operatorname{Sing}_1(m,n)$ for $r\leq d-2$.
The last case follows from Theorem~\ref{thm:singleton} with
$j=d-1$ and $p_j=m-1$.

\medskip
\noindent
\textbf{Outline of the proof.}
By the variational principle of Das, Fishman, Simmons
and Urba\'nski, the dimension problems considered here
reduce to optimizing the time averages of the contraction
rate $\delta(f,t)$ over suitable families of templates.

The main difficulty in proving the upper bounds is that
the local contraction rate may exceed $D_{a,b}$ even
when $f_a(t)<0<f_{b+1}(t)$. To overcome this difficulty,
we construct a continuous piecewise-linear potential
$\Phi$, nonpositive on normalized templates, such that
almost everywhere
\[
 \delta(f,t)-(\Phi\circ f)'(t)
 \leq
 \begin{cases}
 D_{a,b}, & f_a(t)<0<f_{b+1}(t),\\
 mn, & \text{otherwise}.
 \end{cases}
\]
The nonpositivity of the potential allows us to integrate
this inequality and control the average contraction rate
in terms of the proportion of time spent in the window.

The proof proceeds in two steps. In Section~\ref{sec:secants},
we construct a common affine correction that bounds the
contraction rates of all two-block configurations whose
separating indices lie in the prescribed window.
In Section~\ref{sec:potential}, we realize this correction
as the derivative of a piecewise-linear potential and
show that merging adjacent equality blocks among the
negative coordinates, or among the nonnegative coordinates,
cannot decrease the corrected contraction rate.
This reduces the general pointwise estimate to the
two-block case.

Section~\ref{sec:permanent} establishes the dimension formulas
for permanent divergence. Section~\ref{sec:average} proves the
dimension formulas for prescribed divergence proportions, with
the lower bounds obtained by concatenating two-block excursions
of increasing amplitude and intervals on which all coordinates
vanish.

\medskip
\noindent\textit{Acknowledgments.}
The author thanks Tushar Das and David Simmons for useful
discussions and lectures on templates. After the first version
appeared, Tushar Das informed the author that he and David Simmons
had independently established the dimension formula for windows of
length at least two. This communication prompted the author to revisit
and extend the original block-merging argument.

\medskip
\noindent\textit{AI declaration.}
ChatGPT 6 was used as an auxiliary tool to check selected
calculations and mathematical arguments and to improve the
exposition. The author independently reviewed and verified the
mathematical content and takes full responsibility for the manuscript.

\section{Templates and the contraction rate}\label{sec:templates}

Parametric geometry of numbers was introduced by Schmidt and
Summerer \cite{SS} and further developed by Roy \cite{Roy}.
We use the framework of Das, Fishman, Simmons and Urba\'nski,
in which templates model the ordered functions
$h_{A,1},\ldots,h_{A,d}$ up to bounded error
\cite[Theorem~4.2]{DFSU}. The slope restrictions are easiest
to state in terms of partial sums. Wherever the $j$th and
$(j+1)$st coordinates are separated, the sum of the first
$j$ coordinates has one of finitely many allowed slopes.

For $0\leq j\leq d$, define
\[
 Z(j)=\left\{\frac{L_+}{m}-\frac{L_-}{n}:
 L_++L_-=j,\quad
 L_+\in\{0,\ldots,m\},\quad L_-\in\{0,\ldots,n\}\right\}.
\]
An $(m,n)$-template is a continuous, locally piecewise-linear map
$f=(f_1,\ldots,f_d):[0,\infty)\to\R^d$ satisfying the following
conditions:
\begin{enumerate}
 \item $f_1(t)\leq\cdots\leq f_d(t)$ for all $t$;
 \item $-1/n\leq f_i'(t)\leq 1/m$ wherever the derivatives exist;
 \item writing $H_j=f_1+\cdots+f_j$ and $H_0=0$, on every open
 interval on which $f_j<f_{j+1}$, the function $H_j$ is convex
 and its slopes belong to $Z(j)$.
 Here $f_0=-\infty$ and $f_{d+1}=+\infty$.
\end{enumerate}
Since $Z(d)=\{0\}$, the sum $H_d$ is constant. Subtracting $H_d/d$ from each coordinate gives the normalization $H_d=0$ and preserves all slopes and equality blocks. We check bounded-perturbation invariance for the limiting conditions below.

\smallskip
\noindent\textit{Equality blocks and their integer counts.}
A time $t$ is called regular if it belongs to an open interval
on which every coordinate of $f$ is affine and the equality pattern
$f_i=f_{i+1}$ is constant. By local piecewise linearity,
the nonregular times form a locally finite set.

At a regular time, the equality blocks are the maximal consecutive
sets $B=(r,s]\cap\Z$ on which the coordinates of $f$ coincide.
At each block boundary $j$, the partial-sum slope belongs to $Z(j)$.
Consequently there is a unique pair of integers
$L_+(j),L_-(j)$ satisfying
\[
 H_j'=\frac{L_+(j)}{m}-\frac{L_-(j)}{n},
 \qquad L_+(j)+L_-(j)=j.
\]
Indeed, $H_j'= -j/n+dL_+(j)/(mn)$ determines $L_+(j)$, and hence
$L_-(j)$. At the endpoints we have $L_\pm(0)=0$ and
$(L_+(d),L_-(d))=(m,n)$.
For a block $B=(r,s]\cap\Z$, define
\[
 p_B=L_+(s)-L_+(r),\qquad
 q_B=L_-(s)-L_-(r).
\]
Subtracting the count identities $L_+(j)+L_-(j)=j$ at $j=s$ and $j=r$
gives $p_B+q_B=s-r=|B|$. Subtracting the partial-sum slopes
and dividing by $|B|$ gives the common slope of the block,
\begin{equation}\label{eq:slope}
 f_i'=\frac{p_B/m-q_B/n}{p_B+q_B}
     =-\frac1n+\frac{d}{mn}\frac{p_B}{|B|}
 \qquad(i\in B).
\end{equation}
The derivative bounds imply $0\leq p_B\leq |B|$, and hence
$p_B,q_B\geq0$. Moreover, the differences telescope over all equality blocks, giving
\[
 \sum_B p_B=m,\qquad \sum_B q_B=n.
\]

\smallskip
\noindent\textit{The contraction rate.}
Following \cite[Definition~4.5]{DFSU}, within each equality block
$B$ we label the first $p_B$ indices by $+$ and the remaining
$q_B$ indices by $-$. All coordinates in $B$ still have the common slope in~\eqref{eq:slope}. The order of the labels determines which pairs contribute to the contraction rate.

Let $S_+$ and $S_-$ denote the resulting sets of plus and minus
indices, respectively. The local contraction rate at a regular
time is
\[
 \delta(f,t)=\#\{(i_+,i_-)\in S_+\times S_-:i_+<i_-\}.
\]
For the equality blocks $B_1,\ldots,B_\ell$, listed in increasing
order, this gives
\begin{equation}\label{eq:contraction}
 \delta(f,t)=\sum_{r\leq s}p_{B_r}q_{B_s}
            =mn-\sum_{r<s}q_{B_r}p_{B_s}.
\end{equation}
Indeed, a pair is excluded precisely when its minus index lies
in an earlier block than its plus index.
In particular, $0\leq\delta(f,t)\leq mn$.
We set $\delta(f,t)=0$ at nonregular times, which does not affect
the integrals below.

Write $\Delta(f,T)=T^{-1}\int_0^T\delta(f,t)\,dt$ and set
\[
 \underline\delta(f)=\liminf_{T\to\infty}\Delta(f,T),
 \qquad
 \overline\delta(f)=\limsup_{T\to\infty}\Delta(f,T).
\]

\smallskip
\noindent\textit{The variational principle.}
We use the variational principle in the form of
\cite[Theorem~4.6]{DFSU}.
Let $\mathcal B$ be a Borel family of continuous maps
$[0,\infty)\to\R^d$, equipped with the compact-open topology,
which is invariant under bounded perturbations. Then
\begin{equation}\label{eq:VP}
 \begin{aligned}
 \dimH\{A:h_A\in\mathcal B\}
 &=\sup_{f\in\mathcal B\cap\mathcal T_{m,n}}
       \underline\delta(f),\\
 \dimP\{A:h_A\in\mathcal B\}
 &=\sup_{f\in\mathcal B\cap\mathcal T_{m,n}}
       \overline\delta(f),
 \end{aligned}
\end{equation}
where $\mathcal T_{m,n}$ denotes the family of $(m,n)$-templates.

The permanent-divergence conditions and the divergence of a fixed
gap are invariant under bounded perturbations. These conditions
are also Borel. For example, $h_a(t)\to-\infty$ means that for
every integer $N>0$ there exists an integer $T$ such that
\[
 h_a(t)\leq-N
\]
for every rational $t\geq T$. By continuity, the same inequality
then holds for all real $t\geq T$. The same argument applies to
$h_{b+1}(t)\to+\infty$ and to
$h_{i+1}(t)-h_i(t)\to+\infty$.

For the divergence-on-average conditions, suppose
$\sup_{t\geq0}\|h(t)-\widetilde h(t)\|_\infty\leq K$.
For every $R>K$ and $T>0$ we have
\[
 P_{a,b,R+K}(h,T)
 \leq P_{a,b,R}(\widetilde h,T)
 \leq P_{a,b,R-K}(h,T).
\]
Taking the lower limits as $T\to\infty$ and then letting
$R\to\infty$ gives
$P_{a,b}(h)=P_{a,b}(\widetilde h)$.

For fixed $R,T>0$, the map $h\mapsto P_{a,b,R}(h,T)$ is Borel.
Indeed, joint evaluation $(h,t)\mapsto h(t)$ is continuous,
so the indicator of the two threshold inequalities is jointly
Borel, and integration over $[0,T]$ preserves measurability.
Since the integrand is bounded by one, the lower limit as
$T\to\infty$ may be taken along positive integers.
By monotonicity in $R$, the outer limit may also be taken along
positive integers. Consequently, $h\mapsto P_{a,b}(h)$ is Borel.

Thus the families in Theorems~\ref{thm:main}, \ref{thm:singleton}
and~\ref{thm:spectrum}, including their starred versions,
satisfy the hypotheses of \eqref{eq:VP}.
By bounded-perturbation invariance, we may restrict the suprema
in \eqref{eq:VP} to templates normalized by $H_d=0$.

\section{Two-block estimates}\label{sec:secants}

Consider a two-block partition with separating index
$1\leq j\leq d-1$. Write $(p,j-p)$ for the counts of the
left block. Since the total counts are $(m,n)$, the right
block has counts $(m-p,n-j+p)$. By \eqref{eq:contraction},
the contraction rate of this partition is
\begin{equation}\label{eq:delta_j}
 \delta_j(p)=mn-(j-p)(m-p).
\end{equation}
Indeed, the excluded pairs are precisely those whose minus
index lies in the left block and whose plus index lies in
the right block.

Put
\[
 z_j=\frac{mj}{d},\qquad P_j=\lfloor z_j\rfloor,
 \qquad
 \gamma_j=m+j-\lfloor z_j\rfloor-\lceil z_j\rceil.
\]
Let $\nu_L$ and $\nu_R$ be the common coordinate slopes in the left and right blocks. By~\eqref{eq:slope},
\begin{equation}\label{eq:left-slope}
 j\nu_L=\frac{p}{m}-\frac{j-p}{n}=\frac{d}{mn}(p-z_j),
 \qquad (d-j)\nu_R=-j\nu_L.
\end{equation}
The second identity follows because the sum of all coordinate slopes is zero. Thus $p-z_j$ measures the departure of the left block from zero total slope. An admissible plus count is an integer in $[\max\{0,j-n\},\min\{m,j\}]$, and $z_j$ is strictly inside this interval. We use the polynomial~\eqref{eq:delta_j} for all real $p$ because the interpolation below will use its values at integers outside the admissible interval.

\begin{lemma}[Secant majorant at integer points]\label{lem:secant}
For every $1\leq j\leq d-1$ and every $p\in\Z$,
\begin{equation}\label{eq:secant}
 \delta_j(p)-\gamma_j(p-z_j)\leq F_{m,n}(j).
\end{equation}
If $z_j\notin\Z$, equality holds at $p=P_j$ and $p=P_j+1$.
If $z_j\in\Z$, equality holds at $p=P_j$.
\end{lemma}

\begin{proof}
Put $U_j=\lceil z_j\rceil$, and let $S_j$ be the secant line
through $(P_j,\delta_j(P_j))$ and $(U_j,\delta_j(U_j))$ when
$P_j<U_j$, and the tangent line at $P_j$ when $P_j=U_j$.
Since $\delta_j(p)=mn-mj+(m+j)p-p^2$, the difference
$S_j-\delta_j$ is a monic quadratic with roots $P_j,U_j$
(counted with multiplicity). Hence
\[
 S_j(p)-\delta_j(p)=(p-P_j)(p-U_j).
\]
Thus $S_j$ has slope $m+j-P_j-U_j=\gamma_j$. Evaluating at
$p=z_j$ and using
$\{jm/d\}\{jn/d\}=(z_j-P_j)(U_j-z_j)$ gives
$S_j(z_j)=F_{m,n}(j)$. Thus
\begin{equation}\label{eq:secant-difference}
 F_{m,n}(j)+\gamma_j(p-z_j)-\delta_j(p)
 =(p-P_j)(p-U_j).
\end{equation}
Because $P_j$ and $U_j$ are equal or consecutive, the right-hand
side is nonnegative for every integer $p$, with equality exactly
in the stated cases.
\end{proof}

The slope $\gamma_j$ above depends on $j$. We next interpolate the
endpoint slopes to obtain a uniform bound for all intermediate
indices.

\begin{lemma}[Secant interpolation]\label{lem:interp}
Fix integers $1\leq \ell<h\leq d-1$, and let $L(x)$ be
the affine function satisfying
$L(\ell)=\gamma_\ell$ and $L(h)=\gamma_h$.
Then, for every integer $\ell\leq j\leq h$ and every
$p\in\Z$,
\begin{equation}\label{eq:interp}
 \delta_j(p)-L(j)(p-z_j)\leq D_{\ell,h}.
\end{equation}
Moreover, the slope of $L$ has absolute value at most one.
\end{lemma}

\begin{proof}
Write $L(x)=u+vx$. We first show that $|v|\leq1$.
Since
\[
 z_{j+1}-z_j=\frac{m}{d}\in(0,1),
\]
each of $\lfloor z_j\rfloor$ and $\lceil z_j\rceil$
can increase by at most one when $j$ is replaced by $j+1$.
Hence, from
\[
 \gamma_j=m+j-\lfloor z_j\rfloor-\lceil z_j\rceil,
\]
we obtain
\[
 |\gamma_{j+1}-\gamma_j|\leq1.
\]
Therefore
\[
 |\gamma_h-\gamma_\ell|\leq h-\ell,
\]
and since
\[
 v=\frac{\gamma_h-\gamma_\ell}{h-\ell},
\]
we have $|v|\leq1$.

We now prove \eqref{eq:interp}. Fix an integer
$j\in[\ell,h]$ and an integer $p$. Suppose first that $v\geq0$. We keep the plus count $p$
fixed and, for real $x\in[\ell,h]$, consider
\[
 G(x)=\delta_x(p)-L(x)(p-z_x),
\]
where
\[
 \delta_x(p)=mn-(x-p)(m-p),
 \qquad z_x=\frac{mx}{d}.
\]
Thus
\[
 G(x)
 =mn-(x-p)(m-p)
 -(u+vx)\left(p-\frac{mx}{d}\right).
\]
The first two terms are affine in $x$, while the coefficient
of $x^2$ in the last term is $vm/d\geq0$. Hence $G$ is
convex on $[\ell,h]$. Therefore
\[
 G(j)\leq\max\{G(\ell),G(h)\}.
\]
At the endpoints, $L(\ell)=\gamma_\ell$ and
$L(h)=\gamma_h$, so Lemma~\ref{lem:secant} gives
\[
 G(\ell)\leq F_{m,n}(\ell),
 \qquad
 G(h)\leq F_{m,n}(h).
\]
It follows that
\[
 \delta_j(p)-L(j)(p-z_j)
 \leq D_{\ell,h}.
\]

Now suppose that $v<0$. In this case we keep the minus
count
\[
 q=j-p
\]
fixed. Thus, when the separating index is $x$, the
corresponding plus count is $x-q$. Define
\[
 G(x)=\delta_x(x-q)-L(x)(x-q-z_x).
\]
Since
\[
 \delta_x(x-q)=mn-q(m-x+q)
\]
and
\[
 x-q-z_x=\frac{nx}{d}-q,
\]
we have
\[
 G(x)
 =mn-q(m-x+q)
 -(u+vx)\left(\frac{nx}{d}-q\right).
\]
Again, the uncorrected part is affine in $x$, while the
coefficient of $x^2$ is
\[
 -\frac{vn}{d}>0.
\]
Thus $G$ is convex on $[\ell,h]$, and
\[
 G(j)\leq\max\{G(\ell),G(h)\}.
\]
Although these endpoint counts need not be admissible, Lemma~\ref{lem:secant} applies to every integer count. Hence
\[
 G(\ell)\leq F_{m,n}(\ell),
 \qquad
 G(h)\leq F_{m,n}(h),
\]
and therefore
\[
 \delta_j(p)-L(j)(p-z_j)
 \leq D_{\ell,h}.
\]

This proves \eqref{eq:interp}.
\end{proof}

\begin{corollary}\label{cor:endpoint}
For every $1\leq \ell<h\leq d-1$,
\[
 \max_{\ell\leq j\leq h}F_{m,n}(j)
 =\max\{F_{m,n}(\ell),F_{m,n}(h)\}
 =D_{\ell,h},
\]
where the maximum on the left is taken over integer indices.
\end{corollary}

\begin{proof}
Fix $j\in[\ell,h]\cap\Z$ and write $z_j=P_j+\eta$ with
$0\leq\eta<1$. Apply \eqref{eq:interp} at $p=P_j$ and
$p=P_j+1$, multiply the two inequalities by $1-\eta$ and
$\eta$, and add. The correction terms cancel, while
Lemma~\ref{lem:secant} gives
\begin{equation}\label{eq:secant-mean}
 (1-\eta)\delta_j(P_j)+\eta\delta_j(P_j+1)=F_{m,n}(j).
\end{equation}
Hence $F_{m,n}(j)\leq D_{\ell,h}$ for every such $j$.
The reverse inequality is immediate from the two endpoints.
\end{proof}

\smallskip
\noindent\textit{Choice of the affine correction.}
Fix a window $[a,b]$, allowing $a=b$. We seek an affine
correction $C_j$ such that
\[
 Q(j,p):=\delta_j(p)-C_j(p-z_j)\leq D_{a,b}
\]
for $a\leq j\leq b$, while also satisfying $0<C_j<d$ for
every $1\leq j\leq d-1$. To obtain both properties, we enlarge
$[a,b]$ to an interval symmetric about $d/2$.

Set
\[
 r=\min\{a,d-b\}.
\]
Then $[a,b]\subseteq[r,d-r]$, and either $r=a$ or
$d-r=b$. We have the reflection identities
\begin{equation}\label{eq:symmetry}
 F_{m,n}(d-j)=F_{m,n}(j),\qquad
 \gamma_{d-j}=d-\gamma_j.
\end{equation}
Indeed, replacing $j$ by $d-j$ leaves $j(d-j)$ unchanged
and interchanges the two fractional parts in \eqref{eq:F}.
Also, since $z_{d-j}=m-z_j$,
\[
 \lfloor z_{d-j}\rfloor=m-\lceil z_j\rceil,
 \qquad
 \lceil z_{d-j}\rceil=m-\lfloor z_j\rfloor,
\]
which gives the second identity.

There are two cases, according to whether the symmetric interval
$[r,d-r]$ has positive length or consists of the single index
$r=d/2$. Suppose first that $r<d-r$. By
Corollary~\ref{cor:endpoint},
\[
 \max_{r\leq j\leq d-r}F_{m,n}(j)
 =\max\{F_{m,n}(r),F_{m,n}(d-r)\}.
\]
The first identity in \eqref{eq:symmetry} gives
$F_{m,n}(d-r)=F_{m,n}(r)$, and hence
\[
 \max_{r\leq j\leq d-r}F_{m,n}(j)
 =F_{m,n}(r).
\]
Since $a,b\in[r,d-r]$, we obtain
$D_{a,b}\leq F_{m,n}(r)$. Conversely, either $r=a$
or $d-r=b$, so by \eqref{eq:symmetry},
$F_{m,n}(r)$ equals one of
$F_{m,n}(a)$ and $F_{m,n}(b)$. Therefore
\begin{equation}\label{eq:symmetric-window}
 D_{a,b}=F_{m,n}(r).
\end{equation}
If $r=d-r$, then $[a,b]\subseteq[r,r]$, so
$a=b=r$, and \eqref{eq:symmetric-window} is immediate.

For $r<d-r$, let $C_j=c+sj$ be the affine function
interpolating $\gamma_r$ and $\gamma_{d-r}$.
Lemma~\ref{lem:interp} gives $|s|\leq1$. The two interpolation
points are
\[
 (r,\gamma_r)
 \qquad\text{and}\qquad
 (d-r,\gamma_{d-r}).
\]
Setting $j=r$ in the second identity of \eqref{eq:symmetry}
gives $\gamma_{d-r}=d-\gamma_r$. Hence their midpoint is
\[
 \left(
 \frac{r+(d-r)}2,
 \frac{\gamma_r+\gamma_{d-r}}2
 \right)
 =
 \left(\frac d2,\frac d2\right).
\]
Thus the interpolating line passes through $(d/2,d/2)$.

If $r=d-r$, then the second identity in \eqref{eq:symmetry}
gives $\gamma_r=d/2$, and we take $C_j=d/2$. Thus in both
cases
\begin{equation}\label{eq:affine-correction}
 C_j=\frac d2+s\left(j-\frac d2\right),
 \qquad |s|\leq1.
\end{equation}
Therefore, for every $1\leq j\leq d-1$,
\[
 \left|C_j-\frac d2\right|
 \leq\left|j-\frac d2\right|
 <\frac d2,
\]
and hence
\[
 0<C_j<d.
\]

It remains to verify the desired bound for $Q(j,p)$. If
$r<d-r$, Lemma~\ref{lem:interp} applied to $[r,d-r]$
gives
\[
 Q(j,p)\leq D_{r,d-r}
 =F_{m,n}(r)
 =D_{a,b}
\]
for every integer $j\in[r,d-r]$ and every $p\in\Z$.
If $r=d-r$, the same estimate at $j=r$ follows directly
from Lemma~\ref{lem:secant}, since $C_r=\gamma_r$.
Restricting to $j\in[a,b]$, we obtain
\begin{equation}\label{eq:globalC}
 \begin{aligned}
 0<C_j<d
 &\qquad(1\leq j\leq d-1),\\
 Q(j,p)\leq D_{a,b}
 &\qquad(a\leq j\leq b,\ p\in\Z).
 \end{aligned}
\end{equation}

\begin{lemma}[Global corrected-rate bound]\label{lem:global}
For every $1\leq j\leq d-1$ and every admissible plus count
$p$, that is,
\[
 0\leq p\leq m,\qquad 0\leq j-p\leq n,
\]
one has
\[
 Q(j,p)=\delta_j(p)-C_j(p-z_j)\leq mn.
\]
\end{lemma}

\begin{proof}
Put
\[
 q=j-p,\qquad \kappa=\frac{C_j}{d}.
\]
The left block has counts $(p,q)$ and the right block has
counts $(m-p,n-q)$. By \eqref{eq:delta_j},
\[
 mn-\delta_j(p)=q(m-p),
\]
which is exactly the number of pairs with a minus index on the
left and a plus index on the right.

Since $Q(j,p)=\delta_j(p)-C_j(p-z_j)$, we have
\[
 mn-Q(j,p)
 =q(m-p)+\kappa d(p-z_j).
\]
Using $j=p+q$ and $z_j=mj/d$,
\[
 d(p-z_j)
 =np-mq
 =p(n-q)-q(m-p).
\]
Therefore
\[
 \begin{aligned}
 mn-Q(j,p)
 &=q(m-p)
   +\kappa\bigl(p(n-q)-q(m-p)\bigr)\\
 &=(1-\kappa)q(m-p)+\kappa p(n-q).
 \end{aligned}
\]
Since the count $p$ is admissible, we have
\[
 q\geq0,\qquad m-p\geq0,\qquad
 p\geq0,\qquad n-q\geq0.
\]
Moreover, \eqref{eq:globalC} gives $0<C_j<d$, and hence
$0<\kappa<1$. Therefore both terms in
\[
 (1-\kappa)q(m-p)+\kappa p(n-q)
\]
are nonnegative. Thus $mn-Q(j,p)\geq0$, and hence
$Q(j,p)\leq mn$.
\end{proof}

\section{The potential and block-merging argument}
\label{sec:potential}

We seek a potential $\Phi$ whose derivative realizes the
correction $C_j(p-z_j)$ for every two-block configuration. On the hyperplane $\sum_{i=1}^d x_i=0$, define
\begin{equation}\label{eq:potential}
 \Phi(x)=\frac{mn}{d}\sum_{i=1}^d
 \bigl(\alpha_i\min\{x_i,0\}
       +\beta_i\max\{x_i,0\}\bigr),
\end{equation}
where the coefficients $\alpha_i,\beta_i$ will be chosen below.
Each summand vanishes at $x_i=0$, so the two linear pieces agree
there. In particular, $\Phi$ is continuous, piecewise linear, and
Lipschitz.

Consider a normalized two-block template, so that
$f_1+\cdots+f_d=0$, with separating index $j$ and left plus
count $p$. Since the template is normalized and the two block values are
distinct, the coordinates in the left block are negative and
those in the right block are positive. Let $\nu_L,\nu_R$ denote their common
coordinate slopes. By \eqref{eq:left-slope},
\[
 j\nu_L=\frac{d}{mn}(p-z_j),
 \qquad
 (d-j)\nu_R=-j\nu_L.
\]
Define
\[
 \bar\alpha_j=\frac1j\sum_{i=1}^j\alpha_i,
 \qquad
 \bar\beta_j=\frac1{d-j}\sum_{i=j+1}^d\beta_i.
\]
Differentiating \eqref{eq:potential} along this two-block
template gives
\[
 \begin{aligned}
 (\Phi\circ f)'
 &=\frac{mn}{d}
 \left(j\nu_L\bar\alpha_j
 +(d-j)\nu_R\bar\beta_j\right)\\
 &=(p-z_j)(\bar\alpha_j-\bar\beta_j).
 \end{aligned}
\]
Therefore, to obtain
\[
(\Phi\circ f)'=C_j(p-z_j),
\]
it is enough to require
\[
\bar\alpha_j-\bar\beta_j=C_j=c+sj.
\]

We choose two arithmetic progressions. Put
\begin{equation}\label{eq:coeff}
 \begin{gathered}
 u=\max\{-2s,0\},\qquad v=\max\{2s,0\},\\
 \alpha_i=c-u\left(i-\tfrac12\right),\qquad
 \beta_i=v\left(\tfrac{d+1}{2}-i\right).
 \end{gathered}
\end{equation}
A direct calculation gives
\[
 \bar\alpha_j=c-\frac{uj}{2},
 \qquad
 \bar\beta_j=-\frac{vj}{2}.
\]
Since $v-u=2s$, we obtain
\begin{equation}\label{eq:averages}
 \bar\alpha_j-\bar\beta_j
 =c+\frac{v-u}{2}j
 =c+sj
 =C_j.
\end{equation}
Consequently, on every two-block configuration,
\[
 \delta_j(p)-(\Phi\circ f)'
 =\delta_j(p)-C_j(p-z_j)
 =Q(j,p).
\]

We will also use two further properties of these coefficients.
First, both $(\alpha_i)$ and $(\beta_i)$ are nonincreasing,
which will imply $\Phi\leq0$ on the normalized ordered
hyperplane. Second, their common differences are $-u$ and
$-v$. Since $|s|\leq1$, we have $0\leq u,v\leq2$, exactly
the range needed for the block-merging inequality.

\begin{lemma}[Sign of the potential]\label{lem:sign}
If $x_1\leq\cdots\leq x_d$ and $\sum_{i=1}^d x_i=0$,
then $\Phi(x)\leq0$.
\end{lemma}

\begin{proof}
If $x=0$, then $\Phi(x)=0$. Otherwise, let $j$ be the number
of strictly negative coordinates and set
\[
 H=\sum_{i=1}^j x_i<0.
\]
Since the coordinates are ordered and have sum zero,
$1\leq j\leq d-1$ and
\[
 \sum_{i=j+1}^d x_i=-H.
\]

We use the following elementary inequality. If $I$ is a
consecutive set of indices, $(w_i)_{i\in I}$ is nonincreasing,
and $(x_i)_{i\in I}$ is nondecreasing, then, writing
$\bar w=|I|^{-1}\sum_{i\in I}w_i$,
\[
 \sum_{i\in I}w_i x_i-\bar w\sum_{i\in I}x_i
 =
 \frac1{|I|}
 \sum_{\substack{i,k\in I\\ i<k}}
 (w_i-w_k)(x_i-x_k)
 \leq0.
\]
Indeed, for $i<k$ we have $w_i-w_k\geq0$ and
$x_i-x_k\leq0$.

Apply this inequality first to the indices $1,\ldots,j$
with weights $\alpha_i$, and then to the indices
$j+1,\ldots,d$ with weights $\beta_i$. Both coefficient
sequences are nonincreasing by \eqref{eq:coeff}. Hence
\[
 \sum_{i\leq j}\alpha_i x_i\leq\bar\alpha_j H,
 \qquad
 \sum_{i>j}\beta_i x_i\leq-\bar\beta_j H.
\]
Since $x_i<0$ for $i\leq j$ and $x_i\geq0$ for $i>j$,
the definition \eqref{eq:potential} therefore gives
\[
 \begin{aligned}
 \Phi(x)
 &=\frac{mn}{d}
 \left(
 \sum_{i\leq j}\alpha_i x_i
 +\sum_{i>j}\beta_i x_i
 \right)\\
 &\leq
 \frac{mn}{d}
 (\bar\alpha_j-\bar\beta_j)H\\
 &=\frac{mn}{d}C_jH
 <0.
 \end{aligned}
\]
Here we used \eqref{eq:averages} and
$C_j>0$ from \eqref{eq:globalC}, while $H<0$.
\end{proof}

\smallskip
\noindent\textit{The block-merging inequality.}
Let
$\mathcal B=(B_1,\ldots,B_\ell)$ be an ordered partition
of $\{1,\ldots,d\}$ into consecutive blocks. Assign to
each block $B$ nonnegative integer counts $(p_B,q_B)$ such that
\[
 p_B+q_B=|B|,\qquad
 \sum_{B\in\mathcal B}p_B=m,\qquad
 \sum_{B\in\mathcal B}q_B=n.
\]
As in \eqref{eq:slope}, associate to $B$ the common
coordinate slope
\[
 \nu_B
 =-\frac1n+\frac{d}{mn}\frac{p_B}{|B|}.
\]

Fix coefficients $\xi_1,\ldots,\xi_d$, and define their average
on a block $B$ by
\begin{equation}\label{eq:block-average}
 w_B=\frac1{|B|}\sum_{i\in B}\xi_i.
\end{equation}
Define the corrected contraction rate of the partition by
\begin{equation}\label{eq:corrected-rate}
 \mathcal C_\xi(\mathcal B)
 =\sum_{r\leq s}p_{B_r}q_{B_s}
 -\frac{mn}{d}\sum_{B\in\mathcal B}|B|w_B\nu_B.
\end{equation}
Substituting the expression for $\nu_B$ gives
\begin{equation}\label{eq:corrected-counts}
 \mathcal C_\xi(\mathcal B)
 =\sum_{r\leq s}p_{B_r}q_{B_s}
 -\sum_{B\in\mathcal B}p_Bw_B
 +\frac md\sum_{i=1}^d\xi_i.
\end{equation}
The last term is independent of the partition.
Consequently, when two adjacent blocks are merged,
the change in $\mathcal C_\xi$ is determined entirely
by the change in the pair count and in the terms
$p_Bw_B$.

In the merging process, the coefficients $\xi_i$ remain
attached to their original indices, while the counts of
a merged block are obtained by adding the counts of its
constituent blocks.

\begin{lemma}[Block merging]\label{lem:merge}
Let $B_1,B_2$ be two adjacent blocks, listed in increasing
order, with lengths $\ell_i=p_i+q_i>0$.
Suppose that the coefficients on $B_1\cup B_2$ form
an arithmetic progression with common difference
$-\theta$, where $0\leq\theta\leq2$.

Let $\widetilde{\mathcal B}$ be the partition obtained
from $\mathcal B$ by replacing $B_1,B_2$ with their union,
with counts $(p_1+p_2,q_1+q_2)$. Then
\begin{equation}\label{eq:merge}
 \mathcal C_\xi(\widetilde{\mathcal B})
 -\mathcal C_\xi(\mathcal B)
 =
 \left(1-\frac{\theta}{2}\right)p_2q_1
 +\frac{\theta}{2}p_1q_2
 \geq0.
\end{equation}
\end{lemma}

\begin{proof}
We compute the change in the two partition-dependent terms
of \eqref{eq:corrected-counts}.

First consider the pair count. Recall that a plus-minus pair
contributes exactly when the plus label comes before the minus
label. Before the merger, every minus label in $B_1$ comes
before every plus label in $B_2$, so the $p_2q_1$ pairs formed
from these labels are excluded from the count.

After $B_1$ and $B_2$ are replaced by a single block, the
labeling rule places all $p_1+p_2$ plus labels before all
$q_1+q_2$ minus labels. The contribution from within the union
therefore increases by $p_2q_1$. Contributions involving blocks outside
$B_1\cup B_2$ depend only on the total plus and minus counts
of the union, which are preserved by the merger.

Let $w_1,w_2$ denote the block averages defined in
\eqref{eq:block-average} for $B_1,B_2$, and let $w$ denote
the corresponding average on $B_1\cup B_2$.

Since the coefficients on $B_1\cup B_2$ form an arithmetic
progression with common difference $-\theta$, passing from
$B_1$ to the adjacent block $B_2$ lowers the average by
\[
 w_1-w_2=\frac{\theta}{2}(\ell_1+\ell_2).
\]
Indeed, the first term of $B_2$ occurs one step after the last
term of $B_1$, and the mean positions of the two blocks differ
by $(\ell_1+\ell_2)/2$.

Moreover, since $w$ is the average over the union,
\[
 w=\frac{\ell_1w_1+\ell_2w_2}{\ell_1+\ell_2}.
\]
Therefore
\[
 w_1-w
 =\frac{\ell_2}{\ell_1+\ell_2}(w_1-w_2)
 =\frac{\theta\ell_2}{2},
\]
and similarly
\[
 w_2-w
 =-\frac{\ell_1}{\ell_1+\ell_2}(w_1-w_2)
 =-\frac{\theta\ell_1}{2}.
\]

The change in the term $-\sum_B p_Bw_B$ is therefore
\[
 \begin{aligned}
 &p_1w_1+p_2w_2-(p_1+p_2)w\\
 &\qquad=p_1(w_1-w)+p_2(w_2-w)\\
 &\qquad=\frac{\theta}{2}
          (p_1\ell_2-p_2\ell_1)\\
 &\qquad=\frac{\theta}{2}
          (p_1q_2-p_2q_1),
 \end{aligned}
\]
where the last equality follows from
$\ell_i=p_i+q_i$.

Combining this with the increase $p_2q_1$ in the pair count,
we obtain
\[
 \begin{aligned}
 \mathcal C_\xi(\widetilde{\mathcal B})
 -\mathcal C_\xi(\mathcal B)
 &=p_2q_1
   +\frac{\theta}{2}(p_1q_2-p_2q_1)\\
 &=\left(1-\frac{\theta}{2}\right)p_2q_1
   +\frac{\theta}{2}p_1q_2.
 \end{aligned}
\]
Since $0\leq\theta\leq2$ and
$p_1,p_2,q_1,q_2\geq0$, both terms on the right-hand side
are nonnegative. This proves \eqref{eq:merge}.
\end{proof}

The merging operation is purely combinatorial: it acts on the
block partition and its counts at a fixed time; no global template
realizing the merged configurations is required. We now apply
Lemma~\ref{lem:merge} to the equality-block partition at a fixed
regular time, merging all blocks among the negative coordinates
into one block and all blocks among the nonnegative coordinates
into one block.

\begin{lemma}[Pointwise estimate]\label{lem:pointwise}
For every normalized template $f$ and almost every $t$,
\begin{equation}\label{eq:pointwise}
 \delta(f,t)-\frac{d}{dt}\Phi(f(t))
 \leq
 \begin{cases}
 D_{a,b},& f_a(t)<0<f_{b+1}(t),\\
 mn,&\text{otherwise}.
 \end{cases}
\end{equation}
\end{lemma}
\begin{proof}
Since $\Phi$ is Lipschitz and $f$ is locally piecewise linear,
the composition $\Phi\circ f$ is locally absolutely continuous
and hence differentiable almost everywhere. We restrict to
regular times, as defined in Section~\ref{sec:templates}, and
exclude those times $t$ at which, for some $i$,
\[
 f_i(t)=0\qquad\text{and}\qquad f_i'(t)\ne0.
\]
The set of such times has measure zero, since on each affine
piece each coordinate with nonzero slope can vanish at most once.
Thus, at every remaining regular time,
\[
 f_i(t)=0\quad\Longrightarrow\quad f_i'(t)=0
\]
for every $i$.

If $f(t)=0$, the equality-block partition consists of a single
block with counts $(m,n)$. Hence
$\delta(f,t)=mn$ and $(\Phi\circ f)'(t)=0$, so
\eqref{eq:pointwise} holds.

Suppose now that $f(t)\ne0$, and let $j$ be the number of
strictly negative coordinates. Since the coordinates are
ordered and have sum zero,
\[
 1\leq j\leq d-1,\qquad
 f_j(t)<0\leq f_{j+1}(t).
\]
In particular, $j$ is a boundary of the equality-block
partition.

Define
\[
 \xi_i=
 \begin{cases}
 \alpha_i,&i\leq j,\\
 \beta_i,&i>j.
 \end{cases}
\]
Because every zero coordinate has derivative zero at the
chosen time, differentiating \eqref{eq:potential} gives
\begin{equation}\label{eq:potentialderivative}
 (\Phi\circ f)'(t)
 =\frac{mn}{d}\sum_{i=1}^d\xi_i f_i'(t).
\end{equation}
Let $\mathcal B$ be the equality-block partition of $f$ at
time $t$. On each block $B$,
$f_i'(t)=\nu_B$ for every $i\in B$. Thus
\eqref{eq:corrected-rate} and
\eqref{eq:potentialderivative} give
\[
 \delta(f,t)-(\Phi\circ f)'(t)
 =\mathcal C_\xi(\mathcal B).
\]

We now merge all blocks contained in $\{1,\ldots,j\}$ into
one block, and separately merge all blocks contained in
$\{j+1,\ldots,d\}$ into one block. Throughout this process
the coefficients $\xi_i$ remain fixed. On the left they form
an arithmetic progression with common difference $-u$, and
on the right one with common difference $-v$. Since
$u,v\in[0,2]$, Lemma~\ref{lem:merge} shows that each merger
can only increase the corrected contraction rate.

Let $\widetilde{\mathcal B}$ be the resulting two-block
partition, and let $p$ be the plus count of its left block.
Its block counts are
\[
 (p,j-p),\qquad (m-p,n-j+p),
\]
so its contraction rate is $\delta_j(p)$. Let $\nu_L$ and
$\nu_R$ denote the common slopes associated with these two
blocks. By \eqref{eq:left-slope}, their sums of coordinate
slopes are
\[
 \frac{d}{mn}(p-z_j)
 \qquad\text{and}\qquad
 -\frac{d}{mn}(p-z_j),
\]
respectively. By \eqref{eq:corrected-rate}, the term subtracted
from the contraction rate is
\[
 \begin{aligned}
 \frac{mn}{d}
 \left(
 j\bar\alpha_j\nu_L+(d-j)\bar\beta_j\nu_R
 \right)
 &=(p-z_j)(\bar\alpha_j-\bar\beta_j)\\
 &=C_j(p-z_j),
 \end{aligned}
\]
where the last equality follows from \eqref{eq:averages}.
Consequently,
\begin{equation}\label{eq:reduction}
 \delta(f,t)-(\Phi\circ f)'(t)
 =\mathcal C_\xi(\mathcal B)
 \leq\mathcal C_\xi(\widetilde{\mathcal B})
 =Q(j,p).
\end{equation}

The count $p$ is admissible, so Lemma~\ref{lem:global} gives
$Q(j,p)\leq mn$. If in addition
$f_a(t)<0<f_{b+1}(t)$, then the number $j$ of strictly
negative coordinates satisfies $a\leq j\leq b$. In this case
\eqref{eq:globalC} gives
\[
 Q(j,p)\leq D_{a,b}.
\]
Together with \eqref{eq:reduction}, these are exactly the two
bounds in \eqref{eq:pointwise}.
\end{proof}

\section{Dimension bounds for permanent divergence}
\label{sec:permanent}

\begin{proof}[Proof of Theorem~\ref{thm:main}]
\textit{Upper bound.}
Let $f$ be a normalized template satisfying
$f_a(t)\to-\infty$ and $f_{b+1}(t)\to+\infty$.
Choose $t_0$ such that $f_a(t)<0<f_{b+1}(t)$
for all $t\geq t_0$.

Integrating \eqref{eq:pointwise} over $[t_0,T]$
and using Lemma~\ref{lem:sign}, we obtain
\[
 \begin{aligned}
 \int_{t_0}^T\delta(f,t)\,dt
 &\leq D_{a,b}(T-t_0)+\Phi(f(T))-\Phi(f(t_0))\\
 &\leq D_{a,b}(T-t_0)-\Phi(f(t_0)).
 \end{aligned}
\]
Since $0\leq\delta(f,t)\leq mn$, the integral over
$[0,t_0]$ is bounded independently of $T$.
Dividing by $T$ and taking the upper limit as
$T\to\infty$ gives $\overline\delta(f)\leq D_{a,b}$.
By the variational principle \eqref{eq:VP},
\[
 \dimP\mathcal S_{a,b}\leq D_{a,b}.
\]

\smallskip
\noindent\textit{Lower bound.}
For every $a+1\leq k\leq b$, the ordering of the
successive minima gives
$\mathcal S_{k-1,k}\subseteq\mathcal S_{a,b}$,
and hence
$\mathcal S_{k-1,k}^*\subseteq\mathcal S_{a,b}^*$.
The lower bound of
\cite[Theorem~3.14 and Section~3.2.1]{DFSU}
for nontrivially $k$-singular matrices yields
\[
 \begin{aligned}
 \dimH\mathcal S_{a,b}^*
 &\geq
 \max_{a+1\leq k\leq b}
 \max\{F_{m,n}(k-1),F_{m,n}(k)\}\\
 &=\max_{a\leq j\leq b}F_{m,n}(j)
 =D_{a,b},
 \end{aligned}
\]
where the last equality follows from
Corollary~\ref{cor:endpoint}.

Since $\mathcal S_{a,b}^*\subseteq\mathcal S_{a,b}$,
all four dimensions in the theorem lie between
$\dimH\mathcal S_{a,b}^*$ and
$\dimP\mathcal S_{a,b}$.
The preceding bounds therefore prove the theorem.
\end{proof}
\begin{proof}[Proof of Theorem~\ref{thm:singleton}]
Put
\[
 z_j=\frac{mj}{d}.
\]
Recall that
\[
 p_j=\lceil z_j\rceil-1
\]
is the largest integer strictly smaller than $z_j$. Thus
$p_j=P_j=\lfloor z_j\rfloor$ when $z_j\notin\Z$, while
$p_j=P_j-1$ when $z_j\in\Z$.

The upper bound is determined by the number $p$ of plus indices
among the first $j$ positions. We first show that permanent
divergence forces $p\leq p_j$, and then use this to bound the
contraction rate. For the lower bound, we construct a two-block
template attaining the resulting rate.

\smallskip
\noindent\textit{Upper bound.}
Let $f$ be a normalized template satisfying
\[
 f_j(t)\to-\infty,\qquad f_{j+1}(t)\to+\infty.
\]
Then there exists $t_0$ such that
\[
 f_j(t)<f_{j+1}(t)\qquad(t\geq t_0).
\]
By the template axioms, $H_j=f_1+\cdots+f_j$ is convex on
$[t_0,\infty)$ and all of its slopes belong to the finite set
$Z(j)$. Since $H_j$ is convex, its slopes are nondecreasing.
As they take values in a finite set, they can change only
finitely many times. Hence $H_j'$ is constant for all
sufficiently large $t$.

Fix a regular time in this final interval, and let $p$ be the
number of plus indices among $1,\ldots,j$. Since
$f_j(t)<f_{j+1}(t)$, no equality block crosses the index $j$.
For each equality block $B\subseteq\{1,\ldots,j\}$,
\eqref{eq:slope} gives
\[
 \sum_{i\in B}f_i'(t)
 =\frac{p_B}{m}-\frac{q_B}{n}.
\]
Summing over these blocks, whose total plus count is $p$ and
whose total minus count is therefore $j-p$, gives
\[
 H_j'(t)
 =\frac{p}{m}-\frac{j-p}{n}
 =\frac{d}{mn}(p-z_j).
\]
Since $d/(mn)\ne0$, the eventual constant value of $H_j'$
determines a unique integer $p$. Thus the plus count $p$ is
eventually constant.

Moreover,
\[
 H_j(t)=\sum_{i=1}^j f_i(t)
 \leq j f_j(t)\longrightarrow-\infty.
\]
The eventual constant slope of $H_j$ must therefore be negative.
Hence
\[
 p<z_j,
\]
and by the definition of $p_j$,
\[
 p\leq p_j.
\]

Now consider the contraction rate. Among the first $j$ positions
there are $j-p$ minus indices, while among the remaining
positions there are $m-p$ plus indices. Every pair formed by
one of these left minus indices and one of these right plus
indices is excluded from the contraction rate. Since there are
$mn$ plus-minus pairs in total,
\[
 \delta(f,t)\leq mn-(j-p)(m-p).
\]
Since $p\leq p_j$, we have
\[
 j-p\geq j-p_j,\qquad m-p\geq m-p_j.
\]
Therefore, for almost every sufficiently large $t$,
\begin{equation}\label{eq:singleton-upper-rate}
 \delta(f,t)
 \leq mn-(j-p_j)(m-p_j)
 =\delta_j(p_j).
\end{equation}
Consequently,
\[
 \overline\delta(f)
 \leq mn-(j-p_j)(m-p_j).
\]
The variational principle \eqref{eq:VP} now gives
\[
 \dimP\mathcal S_{j,j}
 \leq mn-(j-p_j)(m-p_j),
\]
which is the desired upper bound.

\smallskip
\noindent\textit{Lower bound.}
Since
$\max\{0,j-n\}<z_j<\min\{m,j\}$ and $p_j$ is the largest
integer strictly smaller than $z_j$, the count $p_j$ is admissible.
Set
\[
 w=\frac{p_j}{m}-\frac{j-p_j}{n}
   =\frac{d}{mn}(p_j-z_j)<0,
\]
and define
\[
 f_i(t)=
 \begin{cases}
 wt/j,&1\leq i\leq j,\\
 -wt/(d-j),&j+1\leq i\leq d.
 \end{cases}
\]
For $t>0$ this is a normalized two-block template: the left block
has counts $(p_j,j-p_j)$, the right block has counts $(m-p_j,n-j+p_j)$, and
\eqref{eq:slope} gives their slopes $w/j$ and
$-w/(d-j)$. The only strict gap is after $j$, where
$H_j(t)=wt$ and
$H_j'=p_j/m-(j-p_j)/n\in Z(j)$. Since $w<0$,
$f_j(t)\to-\infty$ and $f_{j+1}(t)\to+\infty$.

By \eqref{eq:delta_j}, its contraction rate is constant:
\[
 \delta(f,t)=\delta_j(p_j)=mn-(j-p_j)(m-p_j)
 \qquad(t>0).
\]
Thus the variational principle gives the matching Hausdorff lower
bound, and together with the packing upper bound proved above and
$\dimH\leq\dimP$ yields \eqref{eq:singleton}.

Finally, every $A\in\mathcal S_{j,j}$ satisfies
$h_{A,j+1}(t)-h_{A,j}(t)\to+\infty$, so it is trivially singular.
Hence $\mathcal S_{j,j}^*=\varnothing$.
\end{proof}

\section{The dimension spectrum for divergence on average}
\label{sec:average}

\begin{proof}[Upper bound in Theorem~\ref{thm:spectrum}]
Let $f$ be a normalized template with $P_{a,b}(f)\geq\rho$,
and fix $R>0$. Put
\[
 E_R(T)=\{t\in[0,T]:f_a(t)\leq-R,\ f_{b+1}(t)\geq R\}.
\]
Then
\[
 |E_R(T)|=T\,P_{a,b,R}(f,T),
\]
and on $E_R(T)$ we have $f_a(t)<0<f_{b+1}(t)$.
Hence Lemma~\ref{lem:pointwise} gives
\[
 \delta(f,t)-(\Phi\circ f)'(t)
 \leq
 \begin{cases}
 D_{a,b},&t\in E_R(T),\\
 mn,&t\in[0,T]\setminus E_R(T),
 \end{cases}
\]
for almost every $t\in[0,T]$. Integrating this inequality over
$[0,T]$ and using the absolute continuity of $\Phi\circ f$, we obtain
\[
 \begin{aligned}
 \int_0^T\delta(f,t)\,dt
 &\leq D_{a,b}|E_R(T)|
      +mn\bigl(T-|E_R(T)|\bigr)\\
 &\qquad+\Phi(f(T))-\Phi(f(0)).
 \end{aligned}
\]
Dividing by $T$ and using Lemma~\ref{lem:sign},
\[
 \Delta(f,T)
 \leq mn-(mn-D_{a,b})P_{a,b,R}(f,T)
      -\frac{\Phi(f(0))}{T}.
\]
Since $mn-D_{a,b}\geq0$, taking the upper limit as
$T\to\infty$ yields
\[
 \overline\delta(f)
 \leq mn-(mn-D_{a,b})
 \liminf_{T\to\infty}P_{a,b,R}(f,T).
\]
Letting $R\to\infty$ and using the definition of
$P_{a,b}(f)$,
\[
 \overline\delta(f)
 \leq mn-(mn-D_{a,b})P_{a,b}(f)
 \leq (1-\rho)mn+\rho D_{a,b}.
\]

Taking the supremum over all such templates in the variational
principle \eqref{eq:VP} gives
\[
 \dimP\mathcal U_{a,b}(\rho)
 \leq(1-\rho)mn+\rho D_{a,b}.
\]
All the other sets appearing in the theorem are subsets of
$\mathcal U_{a,b}(\rho)$, so the same quantity is an upper
bound for all their Hausdorff and packing dimensions.
\end{proof}

\begin{proof}[Lower bound in Theorem~\ref{thm:spectrum}]
By the variational principle \eqref{eq:VP}, it is enough to
construct a normalized template $f$ satisfying
\[
 P_{a,b}(f)=\rho,\qquad
 \underline\delta(f)=(1-\rho)mn+\rho D_{a,b},
\]
and such that no fixed gap $f_{i+1}(t)-f_i(t)$ tends to
infinity. The equality $P_{a,b}(f)=\rho$ places $f$ in the template family
defining $\mathcal L_{a,b}(\rho)$, and the gap condition ensures
that it lies in the starred subfamily
$\mathcal L_{a,b}^*(\rho)$.

For $\rho=0$, take the zero template. Then
$P_{a,b}(f)=0$, all gaps vanish identically, and
$\delta(f,t)=mn$ almost everywhere. Thus assume from now on
that $\rho>0$.

Since
\[
 D_{a,b}=\max\{F_{m,n}(a),F_{m,n}(b)\},
\]
choose an endpoint $j\in\{a,b\}$ such that
\[
 F_{m,n}(j)=D_{a,b}.
\]
Using the notation of Section~\ref{sec:secants}, write
\[
 z_j=\frac{mj}{d}=P_j+\eta,
 \qquad
 P_j=\lfloor z_j\rfloor,
 \qquad
 0\leq\eta<1.
\]

We construct $f$ by concatenating two-block excursions with
intervals on which all coordinates are zero. The $k$th
excursion has length $k^2$, begins and ends at zero, and has
its unique strict internal gap after the fixed index $j$.

\smallskip
\noindent\textit{The excursions.}
For each integer $k\geq3$, construct an excursion
\[
 e^{(k)}=(e_1^{(k)},\ldots,e_d^{(k)})
\]
on the local time interval $0\leq s\leq k^2$. Set
\begin{equation}\label{eq:excursion}
 e_i^{(k)}(s)=
 \begin{cases}
 H_k(s)/j,&1\leq i\leq j,\\
 -H_k(s)/(d-j),&j+1\leq i\leq d,
 \end{cases}
\end{equation}
where $H_k$ will be specified below. By construction,
\[
 H_k(s)=e_1^{(k)}(s)+\cdots+e_j^{(k)}(s).
\]
At an interior regular time, the first $j$ coordinates form
the left equality block. If its plus count is $p$, then
\eqref{eq:left-slope} gives
\begin{equation}\label{eq:excursion-slope}
 H_k'(s)=\frac{d}{mn}(p-z_j).
\end{equation}
For the corresponding two-block partition, the contraction
rate is $\delta_j(p)$ by \eqref{eq:delta_j}.

Suppose first that $0<\eta<1$. Define
\[
 H_k(s)=
 \begin{cases}
 -\dfrac{d\eta}{mn}\,s,
 &0\leq s\leq(1-\eta)k^2,\\[6pt]
 -\dfrac{d(1-\eta)}{mn}(k^2-s),
 &(1-\eta)k^2\leq s\leq k^2.
 \end{cases}
\]
The two formulas agree at $s=(1-\eta)k^2$, where their common
value is $-d\eta(1-\eta)k^2/(mn)$. Hence
\[
 H_k(0)=H_k(k^2)=0,
 \qquad
 H_k(s)<0\quad(0<s<k^2),
\]
and $H_k$ is convex because its slope increases from
$-d\eta/(mn)$ to $d(1-\eta)/(mn)$.

Since $z_j=P_j+\eta$, these slopes correspond by
\eqref{eq:excursion-slope} to plus counts $P_j$ and $P_j+1$,
respectively. The two intervals have lengths $(1-\eta)k^2$
and $\eta k^2$, with contraction rates $\delta_j(P_j)$ and
$\delta_j(P_j+1)$. Thus the mean contraction rate is
\[
 \begin{aligned}
 \frac1{k^2}\int_0^{k^2}\delta(e^{(k)},s)\,ds
 &=(1-\eta)\delta_j(P_j)
   +\eta\delta_j(P_j+1)\\
 &=F_{m,n}(j),
 \end{aligned}
\]
where the last equality follows from \eqref{eq:secant-mean}.

Now suppose that $\eta=0$, so that $z_j=P_j$ is an integer.
Define
\[
 H_k(s)=
 \begin{cases}
 -\dfrac{d}{mn}s,
 &0\leq s\leq k,\\[6pt]
 -\dfrac{dk}{mn},
 &k\leq s\leq k^2-k,\\[6pt]
 -\dfrac{d}{mn}(k^2-s),
 &k^2-k\leq s\leq k^2.
 \end{cases}
\]
The three formulas join continuously and give
\[
 H_k(0)=H_k(k^2)=0,
 \qquad
 H_k(s)<0\quad(0<s<k^2).
\]
The successive slopes $-d/(mn)$, $0$, and $d/(mn)$ are
nondecreasing, so $H_k$ is convex. Since $z_j=P_j$, equation
\eqref{eq:excursion-slope} gives the corresponding plus counts
$P_j-1$, $P_j$, and $P_j+1$ on intervals of lengths $k$,
$k^2-2k$, and $k$, respectively. The associated contraction
rates are $\delta_j(P_j-1)$, $\delta_j(P_j)$, and
$\delta_j(P_j+1)$. Since $\delta_j$ is a quadratic polynomial
with leading coefficient $-1$,
\[
 \delta_j(P_j-1)+\delta_j(P_j+1)
 =2\delta_j(P_j)-2.
\]
Since $z_j=P_j\in\Z$, the equality case of
Lemma~\ref{lem:secant} gives
\[
 \delta_j(P_j)=F_{m,n}(j).
\]
Hence
\[
 \begin{aligned}
 \frac1{k^2}\int_0^{k^2}\delta(e^{(k)},s)\,ds
 &=
 \frac{
 k\delta_j(P_j-1)
 +(k^2-2k)\delta_j(P_j)
 +k\delta_j(P_j+1)}{k^2}\\
 &=F_{m,n}(j)-\frac2k.
 \end{aligned}
\]

We now verify the template axioms for these excursions. Since
$z_j$ lies strictly inside the admissible interval,
\[
 \max\{0,j-n\}<z_j<\min\{m,j\},
\]
the counts $P_j,P_j+1$ in the nonintegral case are admissible.
If $z_j=P_j$ is integral, the same strict inequalities show
that $P_j-1,P_j,P_j+1$ are admissible. Hence the corresponding
coordinate slopes satisfy the derivative bounds.

For $0<s<k^2$, we have $H_k(s)<0$, so
\eqref{eq:excursion} gives
\[
 e_1^{(k)}(s)=\cdots=e_j^{(k)}(s)<0
 <e_{j+1}^{(k)}(s)=\cdots=e_d^{(k)}(s).
\]
Thus, the only strict internal gap is after the index $j$. For the
template condition at this gap, the relevant partial sum is
\[
 H_j(s):=e_1^{(k)}(s)+\cdots+e_j^{(k)}(s)=H_k(s).
\]
If the corresponding plus count is $p$, then
\eqref{eq:excursion-slope} gives
\[
 H_k'(s)
 =\frac{d}{mn}(p-z_j)
 =\frac{p}{m}-\frac{j-p}{n}\in Z(j).
\]
Thus every slope of $H_k$ belongs to $Z(j)$. We have already
seen that these slopes are nondecreasing, so $H_k$ is convex.
Together with the ordering, normalization, and derivative bounds
verified above, this proves that $e^{(k)}$ is an $(m,n)$-template.

Finally,
\[
 e^{(k)}(0)=e^{(k)}(k^2)=0,
\]
so every excursion begins and ends at the zero template.

\smallskip
\noindent\textit{Concatenation.}
For $k\geq3$, set
\[
 T_2=0,\qquad
 T_k=\frac1\rho\sum_{\ell=3}^k\ell^2.
\]
On the $k$th stage $[T_{k-1},T_k]$, place the translated excursion
$e^{(k)}$ on the first $k^2$ units and set all coordinates equal
to zero on the remaining $(\rho^{-1}-1)k^2$ units. Since every
excursion begins and ends at zero, the stages join continuously
and form a normalized $(m,n)$-template. When $\rho=1$ the zero
intervals simply disappear.

Moreover, $f(T_k)=0$ for every $k$. Hence every fixed gap vanishes
at arbitrarily large times and therefore cannot tend to infinity.

\smallskip
\noindent\textit{The time proportion.}
Fix $R>0$. During the $k$th excursion, write
$t=T_{k-1}+s$ with $0\leq s\leq k^2$. Since
$a\leq j\leq b$, \eqref{eq:excursion} gives
\[
 f_a(t)=\frac{H_k(s)}{j},
 \qquad
 f_{b+1}(t)=-\frac{H_k(s)}{d-j}.
\]
Thus the two inequalities
\[
 f_a(t)\leq-R,\qquad f_{b+1}(t)\geq R
\]
hold precisely when
\begin{equation}\label{eq:excursion-threshold}
 H_k(s)\leq-R\max\{j,d-j\}.
\end{equation}

Put
\[
 C_R=R\max\{j,d-j\}.
\]
In both constructions,
\[
 \min_{0\leq s\leq k^2}H_k(s)\longrightarrow-\infty,
\]
so for all sufficiently large $k$ the minimum of $H_k$ is
below $-C_R$. Since $H_k$ decreases linearly from $0$ and later
increases linearly back to $0$, \eqref{eq:excursion-threshold}
fails only on two intervals near the endpoints.

If $0<\eta<1$, the absolute values of the descending and
ascending slopes are $d\eta/(mn)$ and $d(1-\eta)/(mn)$.
Hence the total length of these two endpoint intervals is
\[
 B_R
 :=\frac{mnC_R}{d}
 \left(\frac1\eta+\frac1{1-\eta}\right).
\]
If $\eta=0$, both nonzero endpoint slopes have absolute value
$d/(mn)$, and we set
\[
 B_R:=\frac{2mnC_R}{d}.
\]
Thus $B_R$ is independent of $k$, and for all sufficiently
large $k$ the qualifying time in the $k$th excursion is exactly
$k^2-B_R$. The following zero interval contributes no qualifying
time. Since the $k$th stage has length $k^2/\rho$, its qualifying
proportion is
\[
 \frac{k^2-B_R}{k^2/\rho}
 =\rho\left(1-\frac{B_R}{k^2}\right)
 \longrightarrow\rho.
\]

Let
\[
 T_k=\frac1\rho\sum_{\ell=3}^k\ell^2
\]
be the endpoint of the $k$th stage. The finitely many initial
stages contribute a fixed amount, so there is a constant $A_R$
such that, for all sufficiently large $k$,
\[
 \big|\{t\in[0,T_k]:f_a(t)\leq-R,\ f_{b+1}(t)\geq R\}\big|
 =\sum_{\ell=3}^k\ell^2-kB_R+A_R.
\]
Therefore
\[
 P_{a,b,R}(f,T_k)
 =\rho\left(
 1-\frac{kB_R-A_R}{\displaystyle\sum_{\ell=3}^k\ell^2}
 \right)
 \longrightarrow\rho.
\]

We will use the following observation twice. If
$T_{k-1}\leq T\leq T_k$ and $h=T-T_{k-1}$, then
\[
 \frac{h}{T_{k-1}}
 \leq
 \frac{k^2}{\displaystyle\sum_{\ell=3}^{k-1}\ell^2}
 =O(k^{-1}).
\]
For any integrand $g$ with $0\leq g\leq M$, changing the endpoint
from $T_{k-1}>0$ to $T$ changes its running average by at most
\[
 \frac{Mh}{T_{k-1}}=O(k^{-1}).
\]
Thus convergence along the stage endpoints implies convergence
for arbitrary $T\to\infty$.

Applying this observation to the threshold indicator gives
\[
 \lim_{T\to\infty}P_{a,b,R}(f,T)=\rho
\]
for every fixed $R>0$. Taking the outer limit in the definition
of $P_{a,b}(f)$ gives
\begin{equation}\label{eq:constructed-proportion}
 P_{a,b}(f)=\rho.
\end{equation}

\smallskip
\noindent\textit{The contraction average.}
On every zero interval all coordinates form a single equality
block with counts $(m,n)$. Hence, by
\eqref{eq:contraction},
\[
 \delta(f,t)=mn.
\]

If $0<\eta<1$, the mean contraction rate on every excursion is
\[
 F_{m,n}(j)=D_{a,b}.
\]
Since the excursion occupies the fraction $\rho$ of its stage,
the mean contraction rate of every stage is
\[
 (1-\rho)mn+\rho D_{a,b}.
\]

If $\eta=0$, the mean contraction rate of the $k$th excursion
is
\[
 F_{m,n}(j)-\frac2k
 =D_{a,b}-\frac2k.
\]
Hence the mean contraction rate of the $k$th stage is
\[
 (1-\rho)mn
 +\rho\left(D_{a,b}-\frac2k\right),
\]
which tends to
\[
 (1-\rho)mn+\rho D_{a,b}.
\]

It follows by length-weighted averaging over the completed
stages that
\[
 \Delta(f,T_k)\longrightarrow
 (1-\rho)mn+\rho D_{a,b}.
\]
Applying the observation above to the bounded integrand
$\delta(f,t)$ gives
\[
 \lim_{T\to\infty}\Delta(f,T)
 =(1-\rho)mn+\rho D_{a,b},
\]
and in particular
\begin{equation}\label{eq:constructed-contraction}
 \underline\delta(f)
 =(1-\rho)mn+\rho D_{a,b}.
\end{equation}

Equations~\eqref{eq:constructed-proportion} and
\eqref{eq:constructed-contraction} give the required
properties of $f$. Since no fixed gap tends to infinity,
$f$ belongs to the template family defining
$\mathcal L_{a,b}^*(\rho)$. The variational principle
\eqref{eq:VP} therefore gives
\[
 \dimH\mathcal L_{a,b}^*(\rho)
 \geq(1-\rho)mn+\rho D_{a,b}.
\]

Finally,
\[
 \mathcal L_{a,b}^*(\rho)
 \subseteq\mathcal L_{a,b}(\rho)
 \subseteq\mathcal U_{a,b}(\rho),
 \qquad
 \mathcal L_{a,b}^*(\rho)
 \subseteq\mathcal U_{a,b}^*(\rho)
 \subseteq\mathcal U_{a,b}(\rho).
\]
Together with the packing-dimension upper bound already proved
for $\mathcal U_{a,b}(\rho)$ and the general inequality
$\dimH\leq\dimP$, these inclusions give all of the Hausdorff-
and packing-dimension equalities in
Theorem~\ref{thm:spectrum}.
\end{proof}

\end{document}